\documentclass[a4paper,reqno]{amsart}
\usepackage{geometry}
\usepackage{mathrsfs}
\usepackage{syntonly}
\usepackage{amsmath}
\usepackage{amsthm}
\usepackage{amsfonts}
\usepackage{amssymb}
\usepackage{latexsym}
\usepackage{amscd,amssymb,amsopn,amsmath,amsthm,graphics,amsfonts,mathrsfs,accents,enumerate,verbatim,calc}
\usepackage[dvips]{graphicx}
\usepackage[colorlinks=true,linkcolor=blue,citecolor=blue]{hyperref}
\input xy
\xyoption{all}
\def\del{\delta}
\usepackage[utf8]{inputenc}

\usepackage{color}

\numberwithin{equation}{section}

\def\A{\mathcal{A}}

\def\P{\mathcal{P}}
\def\I{\mathcal{I}}

\def\C{\mathscr{C}}\def\H{\mathcal{H}}
\def\E{\mathbb{E}}
\def\F{\mathbb{F}}
\def\s{\mathfrak{s}}
\def\id{\mathrm{id}}
\def\op{^\mathrm{op}}
\def\Ab{\mathsf{Ab}}
\def\del{\delta}
\def\dr{\ar@{->}[r]}

\def\X{\mathscr{X}}
\def\Y{\mathscr{Y}}

\def\Hom{\mbox{Hom}}

\newcommand{\CC}{{\bf{C}}^{n+2}_{\C}}
\newcommand{\mr}{\hbox{\boldmath$\cdot$}}
\newcommand{\ov}{\overset}
\newcommand{\lra}{\longrightarrow}
\newcommand{\co}{\colon}
\newcommand{\uas}{^{\ast}}            
\newcommand{\sas}{_{\ast}}
\newcommand{\Xd}{\langle X^{\mr},\del\rangle}  
\newcommand{\Yr}{\langle Y^{\mr},\rho\rangle}  

\newcommand{\ush}{^\sharp}           
\newcommand{\ssh}{_\sharp}

\newtheorem{theorem}{Theorem}[section]
\newtheorem{lemma}[theorem]{Lemma}
\newtheorem{corollary}[theorem]{Corollary}

\theoremstyle{definition}
\newtheorem{definition}[theorem]{Definition}

\newtheorem{remark}[theorem]{Remark}
\newtheorem{remark*}[]{Remark}
\newtheorem{example}[theorem]{Example}
\newtheorem{example*}[]{Example}

\newtheorem{condition*}[]{Condition}
\newtheorem{construction}[theorem]{Construction}
\newtheorem{construction*}[]{Construction}

\newtheorem{assumption*}[]{Assumption}

\begin{document}

\title{From right $(n+2)$-angulated categories to $n$-exangulated categories}\footnote{ Panyue Zhou is supported by the National Natural Science Foundation of China (Grant No. 11901190) and by the Scientific Research Fund of Hunan Provincial Education Department (Grant No. 19B239).}
\author{Jian He and Panyue Zhou}
\address{Department of Mathematics, Nanjing University, 210093 Nanjing, Jiangsu, P. R. China}
\email{jianhe30@163.com}
\address{College of Mathematics, Hunan Institute of Science and Technology, 414006 Yueyang, Hunan, China}
\email{panyuezhou@163.com}

\begin{abstract}
The notion of right semi-equivalence in a right $({n+2})$-angulated category is defined in this article. Let $\C$ be an $n$-exangulated category and $\X$ is a strongly covariantly finite subcategory of $\C$. We prove that the standard right $(n+2)$-angulated category $\C/\X$ is right semi-equivalence under a natural assumption.  As an application, we show that a right $(n+2)$-angulated category has an $n$-exangulated structure if and only if the suspension functor is right semi-equivalence. Besides, we also prove that an $n$-exangulated category $\C$ has the structure of a right $(n+2)$-angulated category with right semi-equivalence if and only if for any object $X\in\C$, the morphism $X\to 0$ is a trivial inflation.
\end{abstract}
\keywords{$n$-exangulated categories; extriangulated categories; right $(n+2)$-angulated categories; right triangulated categories; right semi-equivalences}
\subjclass[2020]{18G80; 18G50; 18G25}
\maketitle

\section{Introduction}
Triangulated categories and exact categories are two fundamental structures in algebra, geometry and topology.
They are also important tools in many mathematical branches.
Nakaoka and Palu \cite{NP} recently introduced the notion
of extriangulated categories, whose extriangulated structures are given by $\E$-triangles with some axioms. Triangulated categories and exact categories are extriangulated categories. There are a lot of examples of extriangulated categories which are neither triangulated categories nor exact categories, see \cite{NP,ZZ2,HZZ3,NP1,T}. Beligiannis and Marmaridis \cite{BM} defined the notion of right triangulated category. Informally, a right triangulated category is a triangulated category whose  suspension functor is not necessarily an equivalence. Subsequently, in order to research some homological properties, the notion of a right triangulated category with right semi-equivalence  was introduced by Assem, Beligiannis and Marmaridis \cite{ABM}, these are the right triangulated category whose suspension functor is fully faithful and right dense. We need to pay attention to this relation: a triangulated category is a right triangulated category with right semi-equivalence, but the opposite is not necessarily true, see Example \ref{ex}.
Recently, Tattar \cite{T} proved a right triangulated category with right semi-equivalence is an extriangulated category. More specifically, Tattar \cite{T} gave a characterization of right triangulated categories which are extriangulated as follows.

\begin{theorem}\rm\cite[Theorem A]{T}\label{over}
 Let $(\C, \Sigma, \nabla)$ be a right triangulated category. Then $\C$  has an extriangulated structure if and only if the suspension functor $\Sigma $ is right semi-equivalence.
\end{theorem}
Moreover, Tattar also gave a characterization of extriangulated categories which are right triangulated with right semi-equivalences.
\begin{theorem}\rm\cite[Theorem B]{T}
Let $(\mathscr{C}, \mathbb{E}, \mathfrak{s})$ be an extriangulated category. Then $\mathscr{C}$ has a right triangulated structure and the suspension functor is right semi-equivalence if and only if for any object $ C\in \mathscr{C}$, the morphism $C \rightarrow 0$ is an inflation.

\end{theorem}
In \cite{GKO}, Geiss, Keller and Oppermann introduced
a new type of categories, called $(n+2)$-angulated categories, which generalize triangulated categories: the classical triangulated categories are the special case $n=1$.  Lin \cite{L} introduced right $(n+2)$-angulated categories, the right triangulated categories are the special case $n=1$. Later, Jasso \cite{J} introduced $n$-exact categories which are categories inhabited by certain exact sequences with $n+2$ terms terms, called $n$-exact sequences. The case $n=1$ corresponds to the usual concepts of exact categories. Recently, Herschend, Liu and Nakaoka \cite{HLN} introduced the notion of $n$-exangulated categories. It should be noted that the case $n =1$ corresponds to extriangulated categories. As typical examples we have that $n$-exact and $(n+2)$-angulated categories are $n$-exangulated, see \cite[Proposition 4.5 and Proposition 4.34]{HLN}. However, there are some other examples of $n$-exangulated categories which are neither $(n+2)$-angulated nor $n$-exact, see \cite{HLN,LZ,HZZ2}.

Motivated by the definition of a right triangulated category with
right semi-equivalence and the Tattar's idea, we define right $(n+2)$-angulated categories with right semi-equivalences (see Definition \ref{mdf}), and prove the following.

\begin{theorem}\label{main1} {\rm (}see Theorem \ref{th1} for details{\rm)}
Let $(\C,\E,\s)$ be an $n$-exangulated category and $\X$ be a strongly covariantly finite subcategory of $\C$. Then $\C/\X$ is a right $(n+2)$-angulated category with right semi-equivalence under a reasonable assumption.
\end{theorem}

We also show that  Tattar's result has a higher counterpart:

\begin{theorem}\label{main2} {\rm (}see Corollary \ref{th3} for details{\rm)}
 Let $(\C, \Sigma^n, \Theta)$ be a right $(n+2)$-angulated category. Then $\C$  has the natural structure of an $n$-exangulated category if and only if the $n$-suspension functor $\Sigma^n$ is right semi-equivalence.
\end{theorem}

Moreover, we also gives a way to characterize which $n$-exangulated categories induces
right $(n+2)$-angulated categories with right semi-equivalences.

\begin{theorem}\label{main3} {\rm (}see Theorem \ref{th2} for details{\rm)}
 Let $(\C, \E,\s)$ be an $n$-exangulated category. Then $\C$ induces an
$(n+2)$-angulated structure with right semi-equivalence if and only if for any
object $C\in\C$, the morphism $C\to 0$ is a trivial inflation.
\end{theorem}

This article is organized as follows. In Section 2, we review some elementary definitions and facts on $n$-exangulated categories and right $(n+2)$-angulated categories.
In Section 3, we define a notion of right semi-equivalences for right $(n+2)$-angulated categories,
and we give the proof of Theorem \ref{main1}.
In Section 4, we prove Theorem \ref{main2} and Theorem \ref{main3}.

\section{Preliminaries}
In this section, we briefly review basic concepts and results concerning $n$-exangulated categories and right
$(n+2)$-angulated categories.
\subsection{ $n$-exangulated categories}
{ For any pair of objects $A,C\in\C$, an element $\del\in\E(C,A)$ is called an {\it $\E$-extension} or simply an {\it extension}. We also write such $\del$ as ${}_A\del_C$ when we indicate $A$ and $C$. The zero element ${}_A0_C=0\in\E(C,A)$ is called the {\it split $\E$-extension}. For any pair of $\E$-extensions ${}_A\del_C$ and ${}_{A'}\del{'}_{C'}$, let $\delta\oplus \delta'\in\mathbb{E}(C\oplus C', A\oplus A')$ be the
element corresponding to $(\delta,0,0,{\delta}{'})$ through the natural isomorphism $\mathbb{E}(C\oplus C', A\oplus A')\simeq\mathbb{E}(C, A)\oplus\mathbb{E}(C, A')
\oplus\mathbb{E}(C', A)\oplus\mathbb{E}(C', A')$.

For any $a\in\C(A,A')$ and $c\in\C(C',C)$,  $\E(C,a)(\del)\in\E(C,A')\ \ \text{and}\ \ \E(c,A)(\del)\in\E(C',A)$ are simply denoted by $a_{\ast}\del$ and $c^{\ast}\del$, respectively.

Let ${}_A\del_C$ and ${}_{A'}\del{'}_{C'}$ be any pair of $\E$-extensions. A {\it morphism} $(a,c)\colon\del\to{\delta}{'}$ of extensions is a pair of morphisms $a\in\C(A,A')$ and $c\in\C(C,C')$ in $\C$, satisfying the equality
$a_{\ast}\del=c^{\ast}{\delta}{'}$.}

\begin{definition}\cite[Definition 2.7]{HLN}
Let $\bf{C}_{\C}$ be the category of complexes in $\C$. As its full subcategory, define $\CC$ to be the category of complexes in $\C$ whose components are zero in the degrees outside of $\{0,1,\ldots,n+1\}$. Namely, an object in $\CC$ is a complex $X^{\mr}=\{X_i,d^X_i\}$ of the form
\[ X_0\xrightarrow{d^X_0}X_1\xrightarrow{d^X_1}\cdots\xrightarrow{d^X_{n-1}}X_n\xrightarrow{d^X_n}X_{n+1}. \]
We write a morphism $f^{\mr}\co X^{\mr}\to Y^{\mr}$ simply $f^{\mr}=(f^0,f^1,\ldots,f^{n+1})$, only indicating the terms of degrees $0,\ldots,n+1$.
\end{definition}

\begin{definition}\cite[Definition 2.11]{HLN}
By Yoneda lemma, any extension $\del\in\E(C,A)$ induces natural transformations
\[ \del\ssh\colon\C(-,C)\Rightarrow\E(-,A)\ \ \text{and}\ \ \del\ush\colon\C(A,-)\Rightarrow\E(C,-). \]
For any $X\in\C$, these $(\del\ssh)_X$ and $\del\ush_X$ are given as follows.
\begin{enumerate}
\item[\rm(1)] $(\del\ssh)_X\colon\C(X,C)\to\E(X,A)\ :\ f\mapsto f\uas\del$.
\item[\rm (2)] $\del\ush_X\colon\C(A,X)\to\E(C,X)\ :\ g\mapsto g\sas\delta$.
\end{enumerate}
We simply denote $(\del\ssh)_X(f)$ and $\del\ush_X(g)$ by $\del\ssh(f)$ and $\del\ush(g)$, respectively.
\end{definition}

\begin{definition}\cite[Definition 2.9]{HLN}
 Let $\C,\E,n$ be as before. Define a category $\AE:=\AE^{n+2}_{(\C,\E)}$ as follows.
\begin{enumerate}
\item[\rm(1)]  A pair $\Xd$ is an object of the category $\AE$ with $X^{\mr}\in\CC$
and $\del\in\E(X_{n+1},X_0)$, called an $\E$-attached
complex of length $n+2$, if it satisfies
$$(d_0^X)_{\ast}\del=0~~\textrm{and}~~(d^X_n)^{\ast}\del=0.$$
We also denote it by
$$X_0\xrightarrow{d_0^X}X_1\xrightarrow{d_1^X}\cdots\xrightarrow{d_{n-2}^X}X_{n-1}
\xrightarrow{d_{n-1}^X}X_n\xrightarrow{d_n^X}X_{n+1}\overset{\delta}{\dashrightarrow}.$$
\item[\rm (2)]  For such pairs $\Xd$ and $\langle Y^{\mr},\rho\rangle$,  $f^{\mr}\colon\Xd\to\langle Y^{\mr},\rho\rangle$ is
defined to be a morphism in $\AE$ if it satisfies $(f_0)_{\ast}\del=(f_{n+1})^{\ast}\rho$.

\end{enumerate}
\end{definition}

\begin{definition}\cite[Definition 2.13]{HLN}\label{def1}
 An {\it $n$-exangle} is an object $\Xd$ in $\AE$ that satisfies the listed conditions.
\begin{enumerate}
\item[\rm (1)] The following sequence of functors $\C\op\to\Ab$ is exact.
$$
\C(-,X_0)\xrightarrow{\C(-,\ d^X_0)}\cdots\xrightarrow{\C(-,\ d^X_n)}\C(-,X_{n+1})\xrightarrow{~\del\ssh~}\E(-,X_0)
$$
\item[\rm (2)] The following sequence of functors $\C\to\Ab$ is exact.
$$
\C(X_{n+1},-)\xrightarrow{\C(d^X_n,\ -)}\cdots\xrightarrow{\C(d^X_0,\ -)}\C(X_0,-)\xrightarrow{~\del\ush~}\E(X_{n+1},-)
$$
\end{enumerate}
In particular any $n$-exangle is an object in $\AE$.
A {\it morphism of $n$-exangles} simply means a morphism in $\AE$. Thus $n$-exangles form a full subcategory of $\AE$.
\end{definition}

\begin{definition}\cite[Definition 2.22]{HLN}
Let $\s$ be a correspondence which associates a homotopic equivalence class $\s(\del)=[{}_AX^{\mr}_C]$ to each extension $\del={}_A\del_C$. Such $\s$ is called a {\it realization} of $\E$ if it satisfies the following condition for any $\s(\del)=[X^{\mr}]$ and any $\s(\rho)=[Y^{\mr}]$.
\begin{itemize}
\item[{\rm (R0)}] For any morphism of extensions $(a,c)\co\del\to\rho$, there exists a morphism $f^{\mr}\in\CC(X^{\mr},Y^{\mr})$ of the form $f^{\mr}=(a,f_1,\ldots,f_n,c)$. Such $f^{\mr}$ is called a {\it lift} of $(a,c)$.
\end{itemize}
In such a case, we simple say that \lq\lq$X^{\mr}$ realizes $\del$" whenever they satisfy $\s(\del)=[X^{\mr}]$.

Moreover, a realization $\s$ of $\E$ is said to be {\it exact} if it satisfies the following conditions.
\begin{itemize}
\item[{\rm (R1)}] For any $\s(\del)=[X^{\mr}]$, the pair $\Xd$ is an $n$-exangle.
\item[{\rm (R2)}] For any $A\in\C$, the zero element ${}_A0_0=0\in\E(0,A)$ satisfies
\[ \s({}_A0_0)=[A\ov{\id_A}{\lra}A\to0\to\cdots\to0\to0]. \]
Dually, $\s({}_00_A)=[0\to0\to\cdots\to0\to A\ov{\id_A}{\lra}A]$ holds for any $A\in\C$.
\end{itemize}
Note that the above condition {\rm (R1)} does not depend on representatives of the class $[X^{\mr}]$.
\end{definition}

\begin{definition}\cite[Definition 2.23]{HLN}
Let $\s$ be an exact realization of $\E$.
\begin{enumerate}
\item[\rm (1)] An $n$-exangle $\Xd$ is called an $\s$-{\it distinguished} $n$-exangle if it satisfies $\s(\del)=[X^{\mr}]$. We often simply say {\it distinguished $n$-exangle} when $\s$ is clear from the context.
\item[\rm (2)]  An object $X^{\mr}\in\CC$ is called an {\it $\s$-conflation} or simply a {\it conflation} if it realizes some extension $\del\in\E(X_{n+1},X_0)$.
\item[\rm (3)]  A morphism $f$ in $\C$ is called an {\it $\s$-inflation} or simply an {\it inflation} if it admits some conflation $X^{\mr}\in\CC$ satisfying $d_X^0=f$.
\item[\rm (4)]  A morphism $g$ in $\C$ is called an {\it $\s$-deflation} or simply a {\it deflation} if it admits some conflation $X^{\mr}\in\CC$ satisfying $d_X^n=g$.
\end{enumerate}
\end{definition}

\begin{definition}\cite[Definition 2.27]{HLN}
For a morphism $f^{\mr}\in\CC(X^{\mr},Y^{\mr})$ satisfying $f^0=\id_A$ for some $A=X_0=Y_0$, its {\it mapping cone} $M_f^{\mr}\in\CC$ is defined to be the complex
\[ X_1\xrightarrow{d^{M_f}_0}X_2\oplus Y_1\xrightarrow{d^{M_f}_1}X_3\oplus Y_2\xrightarrow{d^{M_f}_2}\cdots\xrightarrow{d^{M_f}_{n-1}}X_{n+1}\oplus Y_n\xrightarrow{d^{M_f}_n}Y_{n+1} \]
where $d^{M_f}_0=\begin{bmatrix}-d^X_1\\ f_1\end{bmatrix},$
$d^{M_f}_i=\begin{bmatrix}-d^X_{i+1}&0\\ f_{i+1}&d^Y_i\end{bmatrix}\ (1\le i\le n-1),$
$d^{M_f}_n=\begin{bmatrix}f_{n+1}&d^Y_n\end{bmatrix}$.

{\it The mapping cocone} is defined dually, for morphisms $h^{\mr}$ in $\CC$ satisfying $h_{n+1}=\id$.
\end{definition}

\begin{definition}\cite[Definition 2.32]{HLN}
An {\it $n$-exangulated category} is a triplet $(\C,\E,\s)$ of additive category $\C$, additive bifunctor $\E\co\C\op\times\C\to\Ab$, and its exact realization $\s$, satisfying the following conditions.
\begin{itemize}
\item[{\rm (EA1)}] Let $A\ov{f}{\lra}B\ov{g}{\lra}C$ be any sequence of morphisms in $\C$. If both $f$ and $g$ are inflations, then so is $g\circ f$. Dually, if $f$ and $g$ are deflations, then so is $g\circ f$.

\item[{\rm (EA2)}] For $\rho\in\E(D,A)$ and $c\in\C(C,D)$, let ${}_A\langle X^{\mr},c\uas\rho\rangle_C$ and ${}_A\Yr_D$ be distinguished $n$-exangles. Then $(\id_A,c)$ has a {\it good lift} $f^{\mr}$, in the sense that its mapping cone gives a distinguished $n$-exangle $\langle M^{\mr}_f,(d^X_0)\sas\rho\rangle$.
 \item[{\rm (EA2$\op$)}] Dual of {\rm (EA2)}.
\end{itemize}
Note that the case $n=1$, a triplet $(\C,\E,\s)$ is a  $1$-exangulated category if and only if it is an extriangulated category, see \cite[Proposition 4.3]{HLN}.
\end{definition}

\begin{example}
From \cite[Proposition 4.34]{HLN} and \cite[Proposition 4.5]{HLN},  we know that $n$-exact categories and $(n+2)$-angulated categories are $n$-exangulated categories.
There are some other examples of $n$-exangulated categories
 which are neither $n$-exact nor $(n+2)$-angulated, see \cite{HLN,LZ,HZZ2}.
\end{example}

\begin{lemma}\emph{\cite[Proposition 3.6]{HLN}}\label{a2}
\rm Let ${}_A\langle X^{\mr},\delta\rangle_C$ and ${}_B\langle Y^{\mr},\rho\rangle_D$ be distinguished $n$-exangles. Suppose that we are given a commutative square
$$\xymatrix{
 X^0 \ar[r]^{{d_X^0}} \ar@{}[dr]|{\circlearrowright} \ar[d]_{a} & X^1 \ar[d]^{b}\\
 Y^0  \ar[r]_{d_Y^0} &Y^1
}
$$
in $\C$. Then there is a morphism $f^{\mr}\colon\Xd\to\langle Y^{\mr},\rho\rangle$ which satisfies $f^0=a$ and $f^1=b$.
\end{lemma}

\subsection{Right $(n+2)$-angulated categories}
Let $\mathcal{A}$ be an additive category with an endofunctor $\Sigma^n:\mathcal{A}\rightarrow\mathcal{A}$. An $(n+2)$-$\Sigma^n$-$sequence$ in $\mathcal{A}$ is a sequence of morphisms
$$A_0\xrightarrow{f_0}A_1\xrightarrow{f_1}A_2\xrightarrow{f_2}\cdots\xrightarrow{f_{n-1}}A_n\xrightarrow{f_n}A_{n+1}\xrightarrow{f_{n+1}}\Sigma^n A_0.$$
Its {\em left rotation} is the $(n+2)$-$\Sigma^n$-sequence
$$A_1\xrightarrow{f_1}A_2\xrightarrow{f_2}A_3\xrightarrow{f_3}\cdots\xrightarrow{f_{n}}A_{n+1}\xrightarrow{f_{n+1}}\Sigma^n A_0\xrightarrow{(-1)^{n}\Sigma^n f_0}\Sigma^n A_1.$$
A \emph{morphism} of $(n+2)$-$\Sigma^n$-sequences is  a sequence of morphisms $\varphi=(\varphi_0,\varphi_1,\cdots,\varphi_{n+1})$ such that the following diagram commutes
$$\xymatrix{
A_0 \ar[r]^{f_0}\ar[d]^{\varphi_0} & A_1 \ar[r]^{f_1}\ar[d]^{\varphi_1} & A_2 \ar[r]^{f_2}\ar[d]^{\varphi_2} & \cdots \ar[r]^{f_{n}}& A_{n+1} \ar[r]^{f_{n+1}}\ar[d]^{\varphi_{n+1}} & \Sigma^n A_0 \ar[d]^{\Sigma^n \varphi_0}\\
B_0 \ar[r]^{g_0} & B_1 \ar[r]^{g_1} & B_2 \ar[r]^{g_2} & \cdots \ar[r]^{g_{n}}& B_{n+1} \ar[r]^{g_{n+1}}& \Sigma^n B_0
}$$
where each row is an $(n+2)$-$\Sigma^n$-sequence. It is an {\em isomorphism} if $\varphi_0, \varphi_1, \varphi_2, \cdots, \varphi_{n+1}$ are all isomorphisms in $\mathcal{A}$.
\medskip
\begin{definition}\cite[Definition 2.1]{L}\label{d1}
A {\em right} $(n+2)$-\emph{angulated category} is a triple $(\mathcal{A}, \Sigma^n, \Theta)$, where $\mathcal{A}$ is an additive category, $\Sigma^n$ is an endofunctor of $\mathcal{A}$ (which is called the $n$-suspension functor), and $\Theta$ is a class of $(n+2)$-$\Sigma^n$-sequences (whose elements are called right $(n+2)$-angles), which satisfies the following axioms:
\begin{itemize}
\item[(RN1)]
\begin{itemize}
\item[(a)] The class $\Theta$ is closed under isomorphisms, direct sums and direct summands.

\item[(b)] For each object $A\in\mathcal{A}$ the trivial sequence
$$0\rightarrow A\xrightarrow{1_A}A\rightarrow 0\rightarrow\cdots\rightarrow 0\rightarrow 0$$
belongs to $\Theta$.

\item[(c)] Each morphism $f_0\colon A_0\rightarrow A_1$ in $\A$ can be extended to a right $(n+2)$-angle: $$A_0\xrightarrow{f_0}A_1\xrightarrow{f_1}A_2\xrightarrow{f_2}\cdots\xrightarrow{f_{n-1}}A_n\xrightarrow{f_n}A_{n+1}\xrightarrow{f_{n+1}}\Sigma^n A_0.$$
\end{itemize}
\item[(RN2)] If an $(n+2)$-$\Sigma^n$-sequence belongs to $\Theta$, then its left rotation belongs to $\Theta$.

\item[(RN3)] Each solid commutative diagram
$$\xymatrix{
A_0 \ar[r]^{f_0}\ar[d]^{\varphi_0} & A_1 \ar[r]^{f_1}\ar[d]^{\varphi_1} & A_2 \ar[r]^{f_2}\ar@{-->}[d]^{\varphi_2} & \cdots \ar[r]^{f_{n}}& A_{n+1} \ar[r]^{f_{n+1}}\ar@{-->}[d]^{\varphi_{n+1}} & \Sigma^n A_0 \ar[d]^{\Sigma^n \varphi_0}\\
B_0 \ar[r]^{g_0} & B_1 \ar[r]^{g_1} & B_2 \ar[r]^{g_2} & \cdots \ar[r]^{g_{n}}& B_{n+1} \ar[r]^{g_{n+1}}& \Sigma^n B_0
}$$ with rows in $\Theta$ can be completed to a morphism of  $(n+2)$-$\Sigma^n$-sequences.

\item[(RN4)] Given a commutative diagram
$$\xymatrix{
A_0\ar[r]^{f_0}\ar@{=}[d] & A_1 \ar[r]^{f_1}\ar[d]^{\varphi_1} & A_2 \ar[r]^{f_2} & \cdots\ar[r]^{f_{n-1}} & A_{n}\ar[r]^{f_{n}\quad} & A_{n+1} \ar[r]^{f_{n+1}} & \Sigma^n A_0\ar@{=}[d] \\
A_0\ar[r]^{g_0} & B_1 \ar[r]^{g_1}\ar[d]^{h_1} & B_2\ar[r]^{g_2} & \cdots\ar[r]^{g_{n-1}} & B_{n}\ar[r]^{g_{n}\quad} & B_{n+1} \ar[r]^{g_{n+1}} & \Sigma^n A_0\\
& C_2\ar[d]^{h_2} & & & & & \\
& \vdots\ar[d]^{h_{n-1}} & & & & & \\
& C_{n}\ar[d]^{h_{n}} & & & & & \\
& C_{n+1}\ar[d]^{h_{n+1}} & & & & & \\
& \Sigma^n A_1 & & & & & \\
}$$
whose top rows and second column belong to $\Theta$. Then there exist morphisms $\varphi_i\colon A_i\rightarrow B_i\ (i=2,3,\cdots,n+1)$, $\psi_j\colon B_j\rightarrow C_j\ (j=2,3,\cdots,n+1)$ and $\phi_k\colon A_k\rightarrow C_{k-1}\ (k=3,4,\cdots,n+1)$ with the following two properties:

(I) The sequence $(1_{A_1},\varphi_1, \varphi_2,\cdots,\varphi_{n+1})$ is a morphism of $(n+2)$-$\Sigma^n$-sequences.

(II) The $(n+2)$-$\Sigma^n$-sequence
$$A_2\xrightarrow{\left(
                    \begin{smallmatrix}
                      f_2 \\
                      \varphi_2 \\
                    \end{smallmatrix}
                  \right)} A_3\oplus B_2\xrightarrow{\left(
                             \begin{smallmatrix}
                               -f_3 & 0 \\
                               \varphi_3 & -g_2 \\
                               \phi_3 & \psi_2 \\
                             \end{smallmatrix}
                           \right)}
 A_4\oplus B_3\oplus C_2\xrightarrow{\left(
                                       \begin{smallmatrix}
                                         -f_4 & 0 & 0 \\
                                         -\varphi_4 & -g_3 & 0 \\
                                         \phi_4 & \psi_3 & h_2 \\
                                       \end{smallmatrix}
                                     \right)}A_5\oplus B_4\oplus C_3$$
$$\xrightarrow{\left(
                                       \begin{smallmatrix}
                                         -f_5 & 0 & 0 \\
                                         \varphi_5 & -g_4 & 0 \\
                                         \phi_5 & \psi_4 & h_3 \\
                                       \end{smallmatrix}
                                     \right)}\cdots\xrightarrow{\scriptsize\left(\begin{smallmatrix}
             -f_{n} & 0 & 0 \\
             (-1)^{n+1}\varphi_{n} & -g_{n-1} & 0 \\
             \phi_{n} & \psi_{n-1} & h_{n-2} \\
             \end{smallmatrix}
             \right)}A_{n+1}\oplus B_{n}\oplus C_{n-1}$$
$$\xrightarrow{\left(
                                                       \begin{smallmatrix}
                                                         (-1)^{n+1}\varphi_{n+1} &-g_{n} &0 \\
                                                          \phi_{n+1}& \psi_{n}& h_{n-1} \\
                                                       \end{smallmatrix}
                                                     \right)}B_{n+1}\oplus C_{n}\xrightarrow{(\psi_{n+1},\ h_{n})}C_{n+1}\xrightarrow{\Sigma^n f_1\circ h_{n+1}}\Sigma^n A_2 \hspace{10mm}$$
belongs to $\Theta$, and  $h_{n+1}\circ\psi_{n+1}=\Sigma^n f_0\circ g_{n+1}$.
   \end{itemize}
\end{definition}
The notion of a \emph{left $(n+2)$-angulated category} is defined dually.
\vspace{1mm}

If $\Sigma^n$ is an automorphism, it is easy to see that the converse of an axiom (RN2) also holds, thus the right
$(n+2)$-angulated category $(\A,\Sigma^n, \Theta)$ is an $(n+2)$-angulated category in the sense of Geiss-Keller-Oppermann \cite[Definition 1.1]{GKO} and in the sense of Bergh-Thaule \cite[Theorem 4.4]{BT}. If $(\A,\Sigma^n, \Theta)$ is a right $(n+2)$-angulated category, $(\A,\Omega^n, \Phi)$
is a left $(n+2)$-angulated category, $\Omega^n$ is a quasi-inverse of $\Sigma^n$ and $\Theta=\Phi$, then $(\A,\Sigma^n, \Theta)$ is an $(n+2)$-angulated category.

\section{Right $(n+2)$-angulated categories with right semi-equivalences}
We first introduce the concept of  right semi-equivalence in a {\em right} $(n+2)$-\emph{angulated category}.
\begin{definition} \label{mdf}
Let $(\mathcal{C}, \Sigma, \Theta)$ be a {\em right} $(n+2)$-\emph{angulated category}. The endofunctor $\Sigma$ is called a right semi-equivalence if it satisfies the following conditions:
\begin{itemize}
\item[\rm (1)]$\Sigma$ is full;

\item[\rm (2)]$\Sigma$ is faithful;

\item[\rm (3)]$\Sigma$ is right dense, that is, for any morphism $u:U\rightarrow TV$, there exists a right $(n+2)$-angle
$$V\xrightarrow{}U_1\xrightarrow{}U_2\xrightarrow{}\cdots\xrightarrow{}U_n\xrightarrow{}U\xrightarrow{u}\Sigma V.$$ \end{itemize}
\end{definition}

\begin{remark}\label{rk}
Let $(\mathcal{C}, \Sigma, \Theta)$ be a right triangulated category. We recalled the notion of right semi-equivalence was introduced in \cite{ABM}. The suspension functor $\Sigma$ is called right dense if for any right triangle of the form $A\xrightarrow{u}\Sigma B\xrightarrow{v}C^{\prime}\xrightarrow{w}\Sigma A,$
there exists $C\in\mathcal{C}$ such that $C^{\prime}\cong\Sigma C$. The suspension functor $\Sigma$ is called a right semi-equivalence if it is full, faithful and right dense. In fact, if $n=1$, then Definition \ref{mdf} coincides with the definition of right semi-equivalence of right triangulated category in \cite{ABM}.
\proof Let $(\mathcal{C}, \Sigma, \Theta)$ be a right triangulated category and $\Sigma$ be fully faithful. We only need to show that $\Sigma$ is right semi-equivalence if and only if for any morphism $u:U\rightarrow \Sigma V$, there exists a right triangle
$V\xrightarrow{v}W\xrightarrow{w}U \xrightarrow{u} \Sigma V$.

``$\Longrightarrow$" See the proof of Corollary 1.7 in \cite{ABM}.

``$\Longleftarrow$" For any right triangle of the form $U\xrightarrow{u}\Sigma V\xrightarrow{v^{\prime}}W^{\prime}\xrightarrow{w{\prime}}\Sigma U$, there exists a right triangle
$V\xrightarrow{f}M\xrightarrow{g}U \xrightarrow{u} \Sigma V$ by assumption.
 By (RTR2), we have that $U\xrightarrow{u}\Sigma V\xrightarrow{-\Sigma f}\Sigma M \xrightarrow{-\Sigma g} \Sigma U$
is a right triangle. By (RTR4), we have the following commutative diagram
 $$\xymatrix{
U \ar[r]^u \ar@{=}[d] & \Sigma V\ar[r]^{v^{\prime}} \ar@{=}[d] & W^{\prime}\ar[r]^{w^{\prime}} \ar@{-->}[d]^{h}&\Sigma U\ar@{=}[d] &\\
U \ar[r]^{u} & \Sigma V \ar[r]^{-\Sigma f} & \Sigma M\ar[r]^{-\Sigma g} &\Sigma U &}$$
 where the rows are right
triangle in $\mathcal{C}$.
By Corollary 1.5 in \cite{ABM}, we have that $h$ is an isomorphism. This shows that $\Sigma$ is right dense. Moreover, $\Sigma$ is right semi-equivalence.  \qed
\end{remark}

\begin{remark}
Similarly, for a {\em left} $(n+2)$-\emph{angulated category}, we can define left semi-equivalence.
\end{remark}

\begin{example} A $(n+2)$-\emph{angulated category} is a {\em right} $(n+2)$-\emph{angulated category} with right semi-equivalence. But the opposite is not necessarily true. A non-trivial example of a right $(n+2)$-angulated category with right semi-equivalence, see Theorem \ref{th1}.
\end{example}

Now we give some examples of right triangulated categories with right semi-equivalences,
These examples come from \cite[Sections 3.1-3.4]{ABM}, one also can see \cite{Z}.
\begin{example}\label{ex}
\begin{itemize}
\item[\rm (a)]Let $A$ be a finite dimensional $k$-algebra, where $k$ is a commutative ring. We define $$\Lambda=\begin{bmatrix}\begin{matrix} A_0 &&0&\\Q_0 &A_1&&\\&Q_1&A_1 &\\0 &&\ddots&\ddots\end{matrix}\end{bmatrix},$$ where lower triangular matrices have only finitely many non-zero coefficients, $A_i = A, Q_i = DA ={\rm{ Hom}}_k(A, k)$ (with its canonical $A$-$A$-bimodule structure) for all $i\geq0$, addition is the usual addition
of matrices and multiplication is induced from the bimodule structure of $DA$ and the zero morphism
$DA\otimes_{A} DA\rightarrow 0$. Then the stable module category $\rm\overline{mod}\Lambda$ is a right triangulated category, whose
suspension functor is a right semi-equivalence. But the stable module category $\rm\overline{mod}\Lambda$ is not a triangulated category.
\item[\rm (b)] If $A$ is an APR-iterated tilted algebra over an algebraically closed field, then the full subcategory $\H=\H^{-1}(A)$ of $K^{b}(\textrm{proj}A) $ consisting of the complexes having vanishing cohomology in the positive indices is a right triangulated category, whose suspension functor is a right semi-equivalence. But $\H=\H^{-1}(A)$ is not a triangulated category.

\item[\rm (c)] Let $\C$ be a $k$-abelian category, where $k$ is a commutative ring, and $\Y ,\X$ be two k-linear
subcategory of $\C$. Then $\Y $ is said to be $\X$-coresolving \cite{AS} if it satisfies:

(i) $\Y $ contain $\X$;

(ii) If $0\rightarrow A\rightarrow B\rightarrow C\rightarrow 0$ is $\X$-exact, with $A, B \in \Y $, then $C \in \Y $;

(iii) If $0\rightarrow A\rightarrow B\rightarrow C\rightarrow 0$ is $\X$-exact,with $A, C \in \Y $, then $ B\in \Y $.

Assume that $A $ is $\X$-coresolving , each $\X$-monic in $\Y $ is a monomorphism and every object
in $\X$ is $\X$-projective. Then the stable category $\Y/\X$ is a right triangulated category, whose
suspension functor is a right semi-equivalence. But $\Y/\X$ is not a triangulated category.
\end{itemize}
\end{example}
Let $\C$ be an additive category and $\X$ a subcategory of $\C$.
Recall that we say a morphism $f\colon A \to B$ in $\C$ is an $\X$-\emph{monic} if
$$\C(f,X)\colon \C(B,X) \to \C(A,X)$$
is an epimorphism for all $X\in\X$. We say that $f$ is an $\X$-\emph{epic} if
$$\C(X,f)\colon \C(X,A) \to \C(X,B)$$
is an epimorphism for all $X\in\X$.
Similarly,
we say that $f$ is a left $\X$-approximation of $B$ if $f$ is an $\X$-monic and $A\in\X$.
We say that $f$ is a right $\X$-approximation of $A$ if $f$ is an $\X$-epic and $B\in\X$.

We denote by $\C/\X$
the category whose objects are objects of $\C$ and whose morphisms are elements of
$\Hom_{\C}(A,B)/\X(A,B)$ for $A,B\in\C$, where $\X(A,B)$ is the subgroup of $\Hom_{\C}(A,B)$ consisting of morphisms
which factor through an object in $\X$.
Such category is called the (additive) quotient category
of $\C$ by $\X$. For any morphism $f\colon A\to B$ in $\C$, we denote by $\overline{f}$ the image of $f$ under
the natural quotient functor $\C\to\C/\X$.

\begin{definition} \cite[Definition 3.1]{LZ}\label{dd1}
Let $(\C,\E,\s)$ be an $n$-exangulated category. A subcategory $\X$ of $\C$ is called
\emph{strongly contravariantly finite}, if for any object $C\in\C$, there exists a distinguished $n$-exangle
$$B\xrightarrow{}X_1\xrightarrow{}X_2\xrightarrow{}\cdots\xrightarrow{}X_{n-1}\xrightarrow{}X_{n}\xrightarrow{~g~}C\overset{}{\dashrightarrow}$$
where $g$ is a right $\X$-approximation of $C$ and $X_i\in\X$.
Dually we can define \emph{strongly covariantly finite} subcategory.

A strongly contravariantly finite and strongly  covariantly finite subcategory is called \emph{ strongly functorially finite}.
\end{definition}

\begin{construction}\label{con} Assume that $\X$ is a strongly covariantly finite subcategory of an $n$-exangulated category $\C$.

\textbf{Step 1:} For any object $A\in\C$, take a distinguished $n$-exangle $$A\xrightarrow{f}X_1\xrightarrow{}X_{2}\xrightarrow{}\cdots
\xrightarrow{}X_n\xrightarrow{}X_{n+1}\overset{\delta}{\dashrightarrow},$$ where $f$ ia a left $\X$-approximation of $A$ and $X_1,X_2,\cdots,X_{n-1},
X_n\in\X$.

Define $\mathbb{G}(A)=\mathbb{G} A$ to be the image of $X_{n+1}$ in $\C/\X$.

For any morphism $a\in\C(A,A')$, since $f$ is a left $\X$-approximation of $A$ and $X^{\prime}_1\in\X$, we can complete the following commutative diagram:
$$\xymatrix{
A\ar[r]^{f}\ar@{}[dr] \ar[d]^{a} &X_1 \ar[r] \ar@{}[dr]\ar@{-->}[d]^{a_1} &X_2 \ar[r] \ar@{}[dr]\ar@{-->}[d]^{a_2}&\cdot\cdot\cdot \ar[r]\ar@{}[dr] &X_n \ar[r] \ar@{}[dr]\ar@{-->}[d]^{a_n}&X_{n+1} \ar@{}[dr]\ar@{-->}[d]^{b} \ar@{-->}[r]^-{\delta} &\\
{A^{\prime}}\ar[r] &{X^{\prime}_1}\ar[r]&{X^{\prime}_2} \ar[r] &\cdot\cdot\cdot \ar[r] &{X^{\prime} _n}\ar[r]  &{X^{\prime}_{n+1}} \ar@{-->}[r]^-{\delta^\prime} &.}
$$
For any morphism $\bar{a}\in(\C/\X)(A,A')$, define $\mathbb{G}\bar{a}$ to be the image $\bar{b}$ of $b$ in $\C/\X$.

It is not difficult to see that the endofunctor $\mathbb{G}:\C/\X \rightarrow\C/\X$ is well defined.

\textbf{Step 2:} Let $$A_0\xrightarrow{\alpha_0}A_1\xrightarrow{\alpha_1}A_2\xrightarrow{}\cdots\xrightarrow{}A_n\xrightarrow{\alpha_n}A_{n+1}\overset{\eta}{\dashrightarrow}$$ is a distinguished $n$-exangle in $\C$, where $\alpha_0$ is $\X$-monomorphism. It follows that there exists the following commutative diagram of distinguished $n$-exangles
$$\xymatrix{
A_0\ar[r]^{\alpha_0}\ar@{}[dr] \ar@{=}[d] &A_1 \ar[r]^{\alpha_1} \ar@{}[dr]\ar@{-->}[d]^{g_1} &A_2 \ar[r]^{\alpha_2} \ar@{}[dr]\ar@{-->}[d]^{g_2}&\cdot\cdot\cdot \ar[r]\ar@{}[dr] &A_n \ar[r]^{\alpha_n} \ar@{}[dr]\ar@{-->}[d]^{g_{n}}&A_{n+1} \ar@{}[dr]\ar@{-->}[d]^{{\alpha_{n+1}}} \ar@{-->}[r]^-{\eta} &\\
{A_0}\ar[r]^{f_0} &{X_1}\ar[r]^{f_1}&{X_2} \ar[r]^{f_2} &\cdot\cdot\cdot \ar[r] &{X_n}\ar[r]^{f_n}&\mathbb{G} A_0 \ar@{-->}[r]^-{\zeta} &.}
$$

We define $$A_0\xrightarrow{\overline{\alpha_0}} A_1\xrightarrow{\overline{\alpha_1}}A_2\xrightarrow{\overline{\alpha_2}}\cdots\xrightarrow{\overline{\alpha_n}}A_{n+1}\xrightarrow{\overline{\alpha_{n+1}}}\mathbb{G}A_0$$ to be a standard right $(n+2)$-angle in $\C/\X$.

We define
$$\Theta=\{X^\ast|~X^\ast~\textrm{is an}~(n+2){-\mathbb{G}-}\textrm{sequence with}~X^\ast\oplus Y^\ast\cong Z^\ast,~\textrm{where}~Y^\ast~\textrm{is an}~$$
  ~~~~~~~~~~ ~~ ~~~~$(n+2){-\mathbb{G}-}\textrm{sequence and}~Z^\ast~\textrm{is a standard right}~(n+2)\textrm{-angle}\}.$
\end{construction}
\begin{lemma}\label{rig}
Let $\C$ be an $n$-exangulated category and $\X$ is a strongly covariantly finite subcategory of $\C$. Then the quotient category $\C/\X$ is a right $(n+2)$-angulated category with the endofunctor $\mathbb{G}$ and right $(n+2)$-angles defined in Construction \ref{con}.
\end{lemma}
\proof  Since the proof is very similar to \cite[Theorem 3.11]{ZW}, we omit it. Moreover, one also can see \cite[Lemma 3.2]{LZ2}. \qed

\begin{lemma} \rm\cite[Lemma 3.3]{ZW}\label{ZW}
Let $$\xymatrix{
X_0\ar[r]^{f_0}\ar@{}[dr] \ar[d]^{a_0} &X_1 \ar[r]^{f_1} \ar@{}[dr]\ar[d]^{a_1} &X_2 \ar[r]^{f_2} \ar@{}[dr]\ar[d]^{a_2}&\cdot\cdot\cdot \ar[r]\ar@{}[dr] &X_n \ar[r]^{f_n} \ar@{}[dr]\ar[d]^{a_n}&X_{n+1} \ar@{}[dr]\ar[d]^{a_{n+1}} \ar@{-->}[r]^-{\delta} &\\
{Y_0}\ar[r]^{g_0} &{Y_1}\ar[r]^{g_1}&{Y_2} \ar[r]^{g_2} &\cdot\cdot\cdot \ar[r] &{Y _n}\ar[r]^{g_n}  &{Y_{n+1}} \ar@{-->}[r]^-{\eta} &}
$$
be any morphism of distinguished $n$-exangles. Then ${a_0}$ factors through ${f_0}$ if and only if ${a_{n+1}}$ factors through ${g_n}$ if and only if $ (a_0)_{*}{\delta}=(a_{n+1})^{*}{\eta}=0$.

\end{lemma}

\begin{definition}\label{def2}\cite[Definition 3.14 ]{ZW}
Let $(\C,\E,\s)$ be an $n$-exangulated category.
\item[(1)] An object $P\in\C$ is called \emph{projective} if, for any distinguished $n$-exangle
$$A_0\xrightarrow{\alpha_0}A_1\xrightarrow{\alpha_1}A_2\xrightarrow{\alpha_2}\cdots\xrightarrow{\alpha_{n-2}}A_{n-1}
\xrightarrow{\alpha_{n-1}}A_n\xrightarrow{\alpha_n}A_{n+1}\overset{\delta}{\dashrightarrow}$$
and any morphism $c$ in $\C(P,A_{n+1})$, there exists a morphism $b\in\C(P,A_n)$ satisfying $\alpha_n\circ b=c$.
We denote the full subcategory of projective objects in $\C$ by $\P$.
Dually, the full subcategory of injective objects in $\C$ is denoted by $\I$.

\item[(2)] We say that $\C$ {\it has enough  projectives} if
for any object $C\in\C$, there exists a distinguished $n$-exangle
$$B\xrightarrow{\alpha_0}P_1\xrightarrow{\alpha_1}P_2\xrightarrow{\alpha_2}\cdots\xrightarrow{\alpha_{n-2}}P_{n-1}
\xrightarrow{\alpha_{n-1}}P_n\xrightarrow{\alpha_n}C\overset{\delta}{\dashrightarrow}$$
satisfying $P_1,P_2,\cdots,P_n\in\P$. We can define the notion of having \emph{enough injectives} dually.

\end{definition}

\begin{lemma} \rm\cite[Lemma 4.3]{LZ}
Let $(\C,\E,\s)$ be an $n$-exangulated category and $\X$ be a subcategory of $\C$.
For each pair of objects $A$ and $C$ in $\C$, define
$$\F^{\X}(C,A)=
\{A_0\xrightarrow{f}A_1\xrightarrow{}A_2\xrightarrow{}\cdots\xrightarrow{}A_{n-1}
\xrightarrow{}A_n\xrightarrow{}A_{n+1}\overset{\delta}{\dashrightarrow}
\,\mid\, f\ \textrm{is an}\ \X\textrm{-monic}\}.$$
Then $(\C,\F^{\X},\s_{\F^{\X}})$ is an $n$-exangulated category.

Dually, we define for each pair of objects $A$ and $C$ in $\C$
$$\F_{\X}(C,A)=
\{A_0\xrightarrow{}A_1\xrightarrow{}A_2\xrightarrow{}\cdots\xrightarrow{}A_{n-1}
\xrightarrow{}A_n\xrightarrow{g}A_{n+1}\overset{\delta}{\dashrightarrow}
\,\mid\, g\ \textrm{is an}\ \X\textrm{-epic}\}.$$
Then $(\C,\F_{\X},\s_{\F_{\X}})$ is an $n$-exangulated category.

\end{lemma}
We denote the full subcategory of projective objects in  $(\C,\F^{\X},\s_{\F^{\X}})$ by $\rm {Proj}_{\F^{\X}}\C$.

Our first main result is the following.
\begin{theorem}\label{th1}
Let $\C$ be an $n$-exangulated category and $\X$ be a strongly covariantly finite subcategory of $\C$. Then
\begin{itemize}
\item[\rm (1)] The quotient category $\C/\X$ is a right $(n+2)$-angulated category whose suspension functor $\mathbb{G}$ is right dense;
\item[\rm (2)]If $\X\subset \rm {Proj}_{\F^{\X}}\C$, then $\mathbb{G}$ is right semi-equivalence.
\end{itemize}
\end{theorem}

\proof (1) By Lemma \ref{rig}, we know that $\C/\X$ is a right $(n+2)$-angulated category.
Now we show that $\mathbb{G}$ is right semi-equivalence.

For any morphism $\bar{u}: U\rightarrow \mathbb{G}V$ in $\C/\X$, since $\X$ is a strongly covariantly finite subcategory of $\C$, there exists a distinguished $n$-exangle
$$V\xrightarrow{f_0}X_1\xrightarrow{f_1}X_2\xrightarrow{}\cdots\xrightarrow{}X_n\xrightarrow{f_n}\mathbb{G}V\overset{\eta}{\dashrightarrow}$$ where $f_0$ is a left $\X$-approximation of $V$ and $X_1,X_2,\cdots,
X_n\in\X$.
By {\rm (EA2)}, we can observe that $(\id_V,u)$ has a {\it good lift} $f^{\mr}=(\id_V,h_1,\ldots,h_n,u)$, that is, we have the commutative diagram
$$\xymatrix{
V\ar[r]^{g_0}\ar@{}[dr] \ar@{=}[d]^{\id_V} &A_1 \ar[r]^{g_1} \ar@{}[dr]\ar@{-->}[d]^{h_1} &A_2 \ar[r]^{g_2} \ar@{}[dr]\ar@{-->}[d]^{h_2}&\cdot\cdot\cdot \ar[r]\ar@{}[dr] &A_n \ar[r]^{g_n} \ar@{}[dr]\ar@{-->}[d]^{h_{n}}&U \ar@{}[dr]\ar[d]^{{u}} \ar@{-->}[r]^-{u^{\ast}\eta} &\\
{V}\ar[r]^{f_0} &{X_1}\ar[r]^{f_1}&{X_2} \ar[r]^{f_2} &\cdot\cdot\cdot \ar[r] &{X_n}\ar[r]^{f_n}&\mathbb{G} V \ar@{-->}[r]^-{\eta} &}
$$
of distinguished $n$-exangles. We claim that $g_0$ is an $\X$-monic. Indeed, let $a: V\rightarrow X$ be any morphism with $X\in\X$. Since $f_0$ is a left $\X$-approximation, there exists a morphism $b: X_1\rightarrow X$ such that $a=bf_0$ and then $a=bh_1g_0$. This shows that $g_0$ is an $\X$-monic. Thus we have $$V\xrightarrow{\overline{g_0}} A_1\xrightarrow{\overline{g_1}}A_2\xrightarrow{\overline{g_2}}\cdots\xrightarrow{\overline{}}A_{n}\xrightarrow{\overline{g_n}}U\xrightarrow{\overline{u}}\mathbb{G}V$$ is a right $(n+2)$-angle in $\C/\X$ by Construction \ref{con}, hence $\mathbb{G}$ is right dense.

(2) We only need to show that $\mathbb{G}$ is fully faithful by (1).

Let $\bar{u}: U\rightarrow V$ be a morphism in $\C/\X$ with $\mathbb{G}\bar{u}=0$. We have the commutative diagram of distinguished $n$-exangles
$$\xymatrix{
U\ar[r]^{f_u}\ar@{}[dr] \ar[d]^{u} &X_1 \ar[r]^{g_1} \ar@{}[dr]\ar@{-->}[d]^{u_1} &X_2 \ar[r]^{g_2}  \ar@{}[dr]\ar@{-->}[d]^{u_2}&\cdot\cdot\cdot \ar[r]^{g_{n-1}} \ar@{}[dr] &X_n \ar[r]^{g_u}  \ar@{}[dr]\ar@{-->}[d]^{u_n}&\mathbb{G}U \ar@{}[dr]\ar[d]^{\mathbb{G}u=u^{\prime}} \ar@{-->}[r]^-{\delta} &\\
V\ar[r]^{f_v} &{X^{\prime}_1}\ar[r]^{g^{\prime}_1}&{X^{\prime}_2} \ar[r]\ar[r]^{g^{\prime}_2} &\cdot\cdot\cdot \ar[r]^{g^{\prime}_{n-1}} &{X^{\prime} _n}\ar[r]^{g_v}   &\mathbb{G}V \ar@{-->}[r]^-{\delta^\prime} &,}$$
Where ${f_u}$ and ${f_v}$ are left $\X$-approximation, $X_1,X_2,\cdot\cdot\cdot,X_n;~ {X^{\prime}_1},{X^{\prime}_2},\cdot\cdot\cdot,{X^{\prime}_n}\in\X$. Since $\bar{{u}^{\prime}}=\mathbb{G}\bar{u}=0$, there exists $v: X\rightarrow \mathbb{G}V $ and $w: \mathbb{G}U\rightarrow  X $ such that $u^{\prime}=vw$ with $X\in \X$. Since $X\in \rm {Proj}_{\F^{\X}}\C$, there exists $w^{\prime}: X\rightarrow {X^{\prime} _n} $ such that $v=g _vw^{\prime}$, that is, we have the commutative diagram
$$\xymatrix{
 & && &X\ar@{-->}[d]^{w^{\prime}} \ar[dr]^{v} &\mathbb{G}U\ar[l]^{w}\ar@{}[dr]\ar[d]^{\mathbb{G}u=u^{\prime}} &\\
V\ar[r]^{f_v} &{X^{\prime}_1}\ar[r]^{g^{\prime}_1}&{X^{\prime}_2} \ar[r]\ar[r]^{g^{\prime}_2} &\cdot\cdot\cdot \ar[r]^{g^{\prime}_{n-1}} &{X^{\prime} _n}\ar[r]^{g_v}   &\mathbb{G}V \ar@{-->}[r]^-{\delta^\prime} &.}$$
Therefore $u^{\prime}=vw=g _vw^{\prime}w$ and $u^{\prime}$ factors through $g _v$. By Lemma \ref{ZW}, we get $u$ factors through $f _u$, hence $\bar{u}=0$. This shows that $\mathbb{G}$ is faithful.

To prove that $\mathbb{G}$ is full, let $\bar{{u}}:\mathbb{G}U\rightarrow \mathbb{G}V$ be a morphism in $\C/\X$. Since  $\X\subset \rm {Proj}_{\F^{\X}}\C$, We have the commutative diagram of distinguished $n$-exangles by the dual of Lemma \ref{a2}
$$\xymatrix{
U\ar[r]^{f_u}\ar@{}[dr] \ar@{-->}[d]^{u^{\prime}} &X_1 \ar[r]^{g_1} \ar@{}[dr]\ar@{-->}[d]^{u_1} &X_2 \ar[r]^{g_2}  \ar@{}[dr]\ar@{-->}[d]^{u_2}&\cdot\cdot\cdot \ar[r]^{g_{n-1}} \ar@{}[dr] &X_n \ar[r]^{g_u}  \ar@{}[dr]\ar@{-->}[d]^{u_n}&\mathbb{G}U \ar@{}[dr]\ar[d]^{u} \ar@{-->}[r]^-{\delta} &\\
V\ar[r]^{f_v} &{X^{\prime}_1}\ar[r]^{g^{\prime}_1}&{X^{\prime}_2} \ar[r]\ar[r]^{g^{\prime}_2} &\cdot\cdot\cdot \ar[r]^{g^{\prime}_{n-1}} &{X^{\prime} _n}\ar[r]^{g_v}   &\mathbb{G}V \ar@{-->}[r]^-{\delta^\prime} &.}$$
Where ${f_u}$ and ${f_v}$ are left $\X$-approximation, $X_1,X_2,\cdot\cdot\cdot,X_n;~ {X^{\prime}_1},{X^{\prime}_2},\cdot\cdot\cdot,{X^{\prime}_n}\in\X$. Thus $\bar{u}=\mathbb{G}\bar{{u}^{\prime}}$ by Construction \ref{con}.
\qed

As a special case of Theorem \ref{th1} when $n=1$, we have the following.

\begin{corollary}\rm\cite[Lemma 3.3]{T}
Let $\C$ be an extriangulated category and $\X$ be a covariantly finite subcategory of $\C$. If $\X\subset \rm {Proj}_{\F^{\X}}\C$, then the quotient category $\C/\X$ is a right triangulated category with right semi-equivalence.

\end{corollary}

\section{On the relation $n$-exangulated category and right $(n+2)$-angulated category }
\begin{definition}\label{b1} Let $(\mathscr{C},\mathbb{E},\mathfrak{s})$ be an $n$-exangulated category and $X$ be any object in $\C$.
The morphism $X\to 0$ is called \emph{trivial inflation} if it can be embedded in a distinguished $n$-exangle
$$X\rightarrow0\rightarrow0\rightarrow\cdots\rightarrow0
\rightarrow Y\overset{}{\dashrightarrow}.$$
Dually the morphism $0\to X$ is called \emph{trivial deflation} if it can be embedded in a distinguished $n$-exangle
$$Z\rightarrow0\rightarrow0\rightarrow\cdots\rightarrow0
\rightarrow X\overset{}{\dashrightarrow}.$$
\end{definition}

\begin{remark}\label{remark}
In extriangulated categories, note that the inflation and the trivial inflation are the same, and the deflation and the trivial deflation are the same.
\end{remark}

Our second main result is the following.
\begin{theorem}\label{th2}
\rm Let $(\mathscr{C}, \mathbb{E}, \mathfrak{s})$ be an $n$-exangulated category. Then the following statements are equivalent.
 \begin{itemize}
 \item[\rm (1)] There exists a fully faithful additive endofunctor $\Sigma:\C \rightarrow\C$ such that for any morphism $u:U\to \Sigma V$, there exists a distinguished $n$-exangle
$V\rightarrow U_1\rightarrow U_2\rightarrow\cdots\rightarrow U_n
\rightarrow U\overset{}{\dashrightarrow}.$ Moreover, $\mathbb{E}(-,-)\cong\mathscr{C}(-,\Sigma-)$.

\item[\rm (2)] For any object $C \in\C$, the morphism $C \rightarrow 0$ is a trivial inflation.
\item[\rm (3)] There is a right $(n+2)$-angulation of $\C$ which induces the $n$-exangulated category structure $(\mathscr{C}, \mathbb{E}, \mathfrak{s})$.
 \end{itemize}
 \end{theorem}

\proof
$(1)\Rightarrow(2)$ We claim that for any $C\in \mathscr{C}$, the morphism $1_{\Sigma C}\in\mathscr{C}(\Sigma C,\Sigma C)\cong\mathbb{E}(\Sigma C,C)$ is realised by the distinguished $n$-exangle
$$C\rightarrow 0\rightarrow 0\rightarrow\cdots\rightarrow 0
\rightarrow \Sigma C\overset{}{\dashrightarrow}.$$
Indeed, Let $\mathfrak{s}(1_{\Sigma C})=[C\rightarrow E_1\rightarrow E_2\rightarrow\cdots\rightarrow E_n
\rightarrow \Sigma C]$, then we have two exact sequence by Definition \ref{def1}
$$
\C(-,C)\xrightarrow{}\C(-,E_1)\xrightarrow{}\cdots\xrightarrow{}\C(-,\Sigma C)\xrightarrow{~(1_{\Sigma C})\ssh=\rm id~}\E(-,C)\cong\mathscr{C}(-,\Sigma C)
;$$
$$
\C(\Sigma C,-)\xrightarrow{}\C(E_n,-)\xrightarrow{}\cdots\xrightarrow{}\C(C,-)\xrightarrow{~(1_{\Sigma C})\ush=\Sigma~}\E(\Sigma C,-)\cong\mathscr{C}(\Sigma C,\Sigma -)
.$$
Since $\Sigma$ is faithful, then $(1_{\Sigma C})\ssh$ and $(1_{\Sigma C})\ush$ are monomorphic. That is to say $C\rightarrow 0\rightarrow 0\rightarrow\cdots\rightarrow 0
\rightarrow \Sigma C\overset{}{\dashrightarrow}$ is a $n$-exangle. Thus $C\rightarrow E_1\rightarrow E_2\rightarrow\cdots\rightarrow E_n
\rightarrow \Sigma C\overset{}{\dashrightarrow}$ is homotopically equivalent to the object $C\rightarrow 0\rightarrow 0\rightarrow\cdots\rightarrow 0
\rightarrow \Sigma C\overset{}{\dashrightarrow}$ by \cite[Proposition 3.1]{HLN}. This shows  $C \rightarrow 0$ is a trivial inflation.

$(2)\Rightarrow(1)$ For any object $ C\in \mathscr{C}$, since the morphism $C \rightarrow 0$ is a trivial inflation, then we have a distinguished $n$-exangle
$$C\rightarrow0\rightarrow0\rightarrow\cdots\rightarrow0
\rightarrow Y\overset{}{\dashrightarrow}.$$
Notice that $C \rightarrow 0$ is a left $0$-approximation of $C$, then $0$ is a strongly covariantly finite subcategory of an $n$-exangulated category $\C$.
On the other hand, $0\subset \rm {Proj}_{\F^{\X}}\C$. Then $\mathscr{C}\cong\mathscr{C}/0$ has a right $(n+2)$-angulated structure $(\C, \Sigma, \Theta)$ with $\Sigma$ is right semi-equivalence by Construction \ref{con}, Lemma \ref{rig} and Theorem \ref{th1}.

In order to show $\mathbb{E}(-,-)\cong\mathscr{C}(-,\Sigma-)$, we define $$F:\mathscr{C}(C,\Sigma A)\to\mathbb{E}(C,A)$$ $$ f\to f\uas\del$$ for any $A,C\in\C$ and $\del_A\in\mathscr{C}(\Sigma A,A)$. It is not difficult to see that the F is well defined and injective. We only show that F is surjective. Let $\varepsilon\in \mathbb{E}(C,A)$ be realised by $A\rightarrow A_1\rightarrow A_2\rightarrow\cdots\rightarrow A_n
\rightarrow C$, we have the commutative diagram
$$\xymatrix{
A\ar[r]^{}\ar@{}[dr] \ar@{=}[d]^{} &A_1 \ar[r]^{} \ar@{}[dr]\ar@{-->}[d]^{} &A_2 \ar[r]^{} \ar@{}[dr]\ar@{-->}[d]^{}&\cdot\cdot\cdot \ar[r]\ar@{}[dr] &A_n \ar[r]^{} \ar@{}[dr]\ar@{-->}[d]^{}&C \ar@{}[dr]\ar@{-->}[d]^{f} \ar@{-->}[r]^-{f^{\ast}\del_A} &\\
{A}\ar[r]^{} &{0}\ar[r]^{}&{0} \ar[r]^{} &\cdot\cdot\cdot \ar[r] &{0}\ar[r]^{}&\Sigma A \ar@{-->}[r]^-{\del_A} &}
$$
of distinguished $n$-exangles by Lemma \ref{a2} and assumption. Thus we have that $\varepsilon={f^{\ast}\del_A}=:F(f)$.

$(2)\Rightarrow(3)$ Since $\mathbb{E}(-,-)\cong\mathscr{C}(-,\Sigma-)$ by (2), then the right $(n+2)$-angulated structure on $\C$ of Lemma \ref{rig} concides with the $n$-exangulated structure of $\C$. That is to say, the $(n+2)$-angulation of $\C$ (with right semi-equivalence) induces an $n$-exangulation on $\C$.

$(3)\Rightarrow(2)$ Let $(\C, \Sigma, \Theta)$ be a right $(n+2)$-angulation of $\C$. Then for the morphism $1_C \colon C \rightarrow C$ in $\mathscr{C},$ by axioms (RN1) and (RN2), there exists an right $(n+2)$-angle$$C\rightarrow 0\rightarrow\cdots\rightarrow 0\rightarrow \Sigma C\rightarrow\Sigma C$$ in $\Theta$, we have that$$C\rightarrow 0\rightarrow\cdots\rightarrow 0\rightarrow \Sigma C\overset{\delta}{\dashrightarrow}$$ is a distinguished $n$-exangle for some $\delta \in \mathbb{E}(\Sigma C,C)$. So the morphism $C\rightarrow 0$ is a
trivial inflation.

\qed

As a special case of Theorem \ref{th2} when $n=1$, we have the following.
\begin{corollary}\rm\cite[Theorem B]{T}
Let $(\mathscr{C}, \mathbb{E}, \mathfrak{s})$ be an extriangulated category. Then $\mathscr{C}$ has a right triangulated structure $(\C, \Sigma, \Theta)$ with $\Sigma$ is right semi-equivalence if and only if for any object $ C\in \mathscr{C}$, the morphism $C \rightarrow 0$ is an inflation.

\end{corollary}
\proof
This follows from Theorem \ref{th2} and Remark \ref{remark}.  \qed

In general, a right $(n+2)$-angulated category do not have a natural $n$-exangulated structure. As a direct consequence of Theorem \ref{th2}, we have the following.
\begin{corollary}\label{th3}
 Let $(\C, \Sigma^n, \Theta)$ be a right $(n+2)$-angulated category. Then $\C$  has an $n$-exangulated structure $(\mathscr{C}, \mathbb{E}, \mathfrak{s})$ if and only if the suspension functor $\Sigma^n$ is right semi-equivalence.
\end{corollary}

As a special case of Corollary \ref{th3} when $n=1$, we have the following.
\begin{corollary}\rm\cite[Theorem A]{T}\label{over}
 Let $(\C, \Sigma, \Theta)$ be a right triangulated category. Then $\C$  has an extriangulated structure $(\mathscr{C}, \mathbb{E}, \mathfrak{s})$ if and only if the suspension functor $\Sigma $ is right semi-equivalence.
\end{corollary}

\begin{remark}\rm(1) One can get that the Example \ref{ex} can induce extriangulated structures which are not triangulated by Corollary \ref{over}.
\item[\rm (2)]Let $\C$ be an $n$-exangulated category and $\X$ is a strongly covariantly finite subcategory of $\C$. If $\X\subset \rm {Proj}_{\F^{\X}}\C$, then the standard right $(n+2)$-angulated category $\C/\X$ has an $n$-exangulated structure by Theorem \ref{th1} and Corollary \ref{th3}.

\end{remark}

\end{document}